\documentclass[12pt]{article}
\usepackage{amssymb}
\usepackage{amsmath}
\usepackage{amsthm}
\usepackage{anysize}\marginsize{45mm}{45mm}{40mm}{50mm}
\usepackage{bbm}
\usepackage{caption}
\usepackage[usenames,dvipsnames]{xcolor}
\usepackage{cite}
\usepackage{float}
\usepackage{graphicx}
\usepackage{indentfirst}
\usepackage{latexsym}
\usepackage{latexsym, bm}
\usepackage{mathrsfs}
\usepackage{subfigure}

\usepackage{enumerate}

\newtheoremstyle{lemma}{\topsep}{\topsep}%
     {}
     {}
     {\bfseries}
     {}
     {0.1em}
     {\thmname{#1}\thmnumber{ #2}\thmnote{ #3}}
\theoremstyle{lemma}  

\newtheorem{theorem}{Theorem}[section]    
\newtheorem{lemma}[theorem]{Lemma}
\newtheorem{corollary}[theorem]{Corollary}
\newtheorem{proposition}[theorem]{Proposition}

\newtheorem{example}[theorem]{Example}

\numberwithin{equation}{section}

\title{On anti-Kekul\'{e} and $s$-restricted matching preclusion problems}

\author{Huazhong L\"{u}$^{1,2}$\thanks{Corresponding author.}, Xianyue Li$^{1}$ and Heping Zhang$^{1}$ \\
{\footnotesize 1. School of Mathematics and Statistics, Lanzhou University, }\\
{\footnotesize Lanzhou, Gansu 730000, P.R. China}\\
{\footnotesize 2. School of Mathematical Sciences, University of Electronic Science and Technology of China,} \\ {\footnotesize Chengdu, Sichuan 610054, P.R. China}\\
{\footnotesize E-mails: lvhz@uestc.edu.cn, lixianyue@lzu.edu.cn, zhanghp@lzu.edu.cn}}
\date{}
\begin{document}

\maketitle
\begin{abstract}

The anti-Kekul\'{e} number of a connected graph $G$ is the smallest number of edges whose deletion results in a
connected subgraph having no Kekul\'{e} structures (perfect matchings). As a common generalization of (conditional) matching preclusion number and anti-Kekul\'{e} number of a graph $G$, we introduce $s$-restricted matching preclusion number of $G$ as the smallest number of edges whose deletion results in a subgraph without perfect matchings such that each component has at least $s+1$ vertices. In this paper, we first show that conditional matching preclusion problem and anti-Kekul\'{e} problem are NP-complete, respectively, then generalize this result to $s$-restricted matching preclusion problem. Moreover, we give some sufficient conditions to compute $s$-restricted matching preclusion numbers of regular graphs. As applications, $s$-restricted matching preclusion numbers of complete graphs, hypercubes and hyper Petersen networks are determined.
%

\vskip 0.05 in

\end{abstract}

\section{Introduction}

Let $G=(V,E)$ be a simple and connected graph. Let $N(v)$ be the set of neighbors of a vertex $v$ and $d(v)=\vert N(v)\vert$, the {\em degree} of $v$. The minimum degree of $G$ is denoted by $\delta(G)$. A {\em
matching} $M$ of $G$ is a set of pairwise nonadjacent edges of $G$.
The vertices of $G$ incident to the edges of $M$ are called {\em
saturated} by $M$; the others are {\em unsaturated}. A matching in
$G$ of maximum cardinality is called a {\em maximum matching}. The
cardinality of a maximum matching is called the {\em matching
number} of $G$, denoted by $\nu(G)$. A {\em perfect matching} is a
matching of cardinality $\vert V(G)\vert/2$. An {\em almost perfect
matching} is a matching covering all but one vertex of $G$. Let $F$
be a set of edges of $G$. If $G-F$ has neither a perfect matching nor
an almost perfect matching, then we call $F$ a {\em matching preclusion
set} of $G$. The {\em matching preclusion number} of $G$, denoted by
$mp(G)$, is the minimum cardinality over all matching preclusion sets of
$G$. A matching preclusion set of a graph $G$ with even order is {\em trivial} if all its edges are incident to a vertex of $G$. Based on the
definition, we set $mp(G)=0$ if $G$ has neither a perfect matching nor
an almost perfect matching. For other standard graph notations and
terminologies not defined here please refer to \cite{Bondy}.

In organic molecule graphs, perfect matchings correspond to
Kekul\'{e} structures, which play an important role in analyzing
resonant energy and stability of hydrocarbon compounds. In
\cite{Vukicevic1}, Vuki\v{c}evi\'{c} and Trinajsti\'{c} proposed the {\em anti-Kekul\'{e} number} of a connected graph $G$, denoted by $ak(G)$, as the smallest number of edges such that after deleting these edges of $G$ the
resulting graph remains connected but has no Kekul\'{e}
structures (perfect matchings). For convenience, we call such a set of edges of $G$ an {\em anti-Kekul\'{e} set}. Anti-Kekul\'{e} numbers of some chemical graphs were studied, such as a hexagonal system \cite{Cai}, cata-condensed benzenoids \cite{Vukicevic2}, fullerene graphs \cite{Yang} and $(4,5,6)$-fullerenes \cite{Zhao}. It is noticeable that not every graph necessarily has an anti-Kekul\'{e} set, such as $K_2$ and even cycles. Recently, graphs which do not have anti-Kekul\'{e} set were characterized and constructed by Wu and Zhang \cite{Wu}.

The concept of matching preclusion was first introduced by Brigham
et al.\cite{Brigham} as a measure of robustness of interconnection
networks under the condition of edge failure. In the same
paper, the authors showed that it will be more robust under edge failure if each vertex has a special matching vertex at any time, and they also
determined the matching preclusion number of Petersen graph,
complete graphs $K_n$, complete bipartite graphs $K_{n,n}$ and hypercubes. Recently, matching preclusion numbers for Cayley graphs generated by transposition trees and ($n,k$)-star graphs \cite{Eddie2}, tori (including related Cartesian product graphs) \cite{Eddie3}, binary de Bruijn graphs \cite{Lin}, $n$-grid graphs \cite{Ding} and data center networks \cite{Lu2} have been determined.

In large networks failure is inevitable, but it is unlikely that
all the edges incident to a common vertex are all faulty
simultaneously. Thus, it is meaningful to consider
matching preclusion of a graph with some restriction on the order of components after edge deletion. Motivated by this, Cheng et al. \cite{Eddie1} considered
conditional matching preclusion set (resp. number) of a graph $G$.
The {\em conditional matching preclusion number} of a graph $G$,
denoted by $mp_{1}(G)$, is the minimum number of edges whose
deletion leaves the resulting graph with no isolated vertices and
without a perfect matching or an almost perfect matching. This problem has been solved for complete graphs, complete bipartite graphs and hypercubes \cite{Eddie1}, Cayley graphs generated by 2-trees and hyper Petersen networks \cite{Eddie5}, HL-graphs \cite{Son}, $k$-ary $n$-cubes \cite{Wang}, balanced hypercubes \cite{Lu} and cube-connected cycles \cite{Li}.

From definitions of the anti-Kekul\'{e} number, the matching preclusion number, and the conditional matching preclusion number, it can be seen that the common point is that edge deletion results in the remaining graph no longer possessing a perfect matching, and the difference is variable requirements about the orders of components of the remaining graph. Motivated by $k$-restricted edge-connectivity \cite{Fabrega}, it is natural to ask how the minimum size of matching preclusion sets $F$ changes if each component of $G-F$ contains at least $s+1$ vertices, where $s$ is a nonnegative integer. This is closely related to the changing and unchanging of invariants studied in many areas \cite{Harary}. Based on this fact, we generalize them as follows.

Let $s$ be a nonnegative integer and $F$ an edge subset of $G$. If $G-F$ has neither a perfect matching nor an almost perfect matching, and each component of $G-F$ has at least $s+1$ vertices, then $F$ is called an {\em $s$-restricted matching preclusion set} of $G$. The $s$-{\em restricted matching preclusion number}, denoted by $mp_{s}(G)$, is the minimum cardinality over all $s$-restricted matching preclusion sets of $G$. It is suitable to make a convention that $s\leq|V(G)|-1$ for a given graph $G$.

Based on the definition, we set $mp_{s}(G)=0$ if $G$ contains no components of order at most $s$ ($s>0$), and has neither perfect matchings nor almost perfect matchings. We define $mp_{s}(G)=+\infty$ if the $s$-restricted matching preclusion set does not exist, that is, we can not delete edges to satisfy the conditions in the definition.

It is noticeable that $0$-restricted matching preclusion problem is equivalent to matching preclusion problem, $1$-restricted matching preclusion problem is the conditional matching preclusion problem and $(|V(G)|-1)$-restricted matching preclusion problem is clearly the anti-Kekul\'{e} problem. Thus, for notation consistency, we still use $mp(G)$ to denote $mp_{0}(G)$.

In \cite{Lacroix}, Lacroix et al. considered matching preclusion problem of graphs with a perfect matching and proved its NP-completeness for bipartite graphs. In view of the similarity of the anti-Kekul\'{e} problem, the conditional matching preclusion problem and the $s$-restricted matching preclusion problem, it is natural to ask what are the complexities of these problems for general graphs. In this paper, we solve this question by showing that they are all NP-complete for bipartite graphs. Additionally, we give some methods to compute $s$-restricted matching preclusion numbers for regular graphs, and as applications, give examples of calculating $s$-restricted matching preclusion number of three kinds of interconnection networks.

The rest of the paper is organized as follows. In Section 2, we prove NP-completeness of the conditional matching preclusion problem by reducing MBPMP to it, and then, as a direct corollary, obtain NP-completeness of the anti-Kekul\'{e} problem. Moreover, we obtain relationships of $s$-restricted matching preclusion numbers concerning anti-Kekul\'{e} number when $s$ increases, and also obtain NP-completeness of the $s$-restricted matching preclusion problem. In Section 3, we present some sufficient conditions to calculate $s$-restricted ($s\geq2$) matching preclusion number of regular graphs. Applications of determining $s$-restricted matching preclusion numbers of complete graphs, hypercubes and hyper Petersen networks are given in Section 4. Finally, we conclude this paper in Section 5.

\section{Complexity results}

Let $G=(V,E)$ be a graph with matching number $\nu(G)$. A {\em
blocker} for $G$ is a set of edges whose removal results in
matching number of $G$ smaller than $\nu(G)$. The decision problem of {\em
blocker problem} (BP) is defined as follows. Given $G$ and a positive integer $k$, does there exist an edge
subset $B$ of $E$ with $\vert B\vert\leq k$ such that $B$ is a
blocker for $G$? In \cite{Zenklusen}, Zenklusen et al. showed that BP of a bipartite graph is NP-complete. Later, Lacroix et al. \cite{Lacroix} considered a special case of BP, called {\em minimum blocker perfect
matching problem} (MBPMP), where $G$ is a graph with a perfect matching. They
proved the following statement.

\begin{lemma}\cite{Lacroix}{\bf .} MBPMP of a bipartite graph is NP-complete.
\end{lemma}

The decision problem of matching preclusion problem is defined as follows.

\textsc{Matching preclusion problem:}

\vskip 0.00 in

\textit{Instance:} A nonempty connected graph $G=(V,E)$ and a positive integer $k$.

\vskip 0.00 in

\textit{Question:} Does there exist a subset $B\subseteq E$ with $\vert
B\vert\leq k$ such that $G-B$ has neither perfect matchings nor almost perfect matchings?

Clearly, MBPMP is equivalent to matching preclusion problem when we restrict our consideration to graphs with a perfect matching. Thus, the following result is straightforward.

\begin{corollary}{\bf .} Matching preclusion problem of a bipartite graph is NP-complete.
\end{corollary}

We shall prove the NP-completeness of conditional matching preclusion problem by reducing from MBPMP to it. We present its decision problem as follows.

\vskip 0.1 in

\textsc{Conditional matching preclusion problem:}

\vskip 0.00 in

\textit{Instance:} A nonempty connected graph $G=(V,E)$ having a perfect matching and a positive integer $k$.

\vskip 0.00 in

\textit{Question:} Does there exist a subset $B\subseteq E$ with $\vert
B\vert\leq k$ such that $G-B$ has neither isolated vertices nor perfect matchings?

\begin{theorem}\label{NPC-conditional}{\bf .}
Conditional matching preclusion problem of a bipartite graph is NP-complete.
\end{theorem}
\begin{proof} Obviously, conditional matching preclusion problem is in NP, because we can check in polynomial time whether a set of edges is a conditional matching preclusion set. We will prove NP-hardness of conditional matching preclusion problem by reducing MBPMP to it in polynomial time.

Let $G=(U\cup V,E)$ be a bipartite graph with bipartition $U$ and
$V$ such that $\vert U\vert=\vert V\vert=t$. Suppose that $G$ has a
perfect matching. Let $u_{1},u_{2},\ldots,u_{t}$ (resp.
$v_{1},v_{2},\ldots,v_{t}$) denote the vertices in $U$ (resp.
$V$). The graph $G'=(U'\cup V',E')$ is constructed from $G$ as
follows (see Fig. \ref{g1}). $U'=U\cup \{u',u''\}$, $V'=V\cup \{v',v''\}$, where
$u',u'',v'$ and $v''$ are new added vertices. $E'=E\cup \{u'v:v\in
V\}\cup \{uv':u\in U\}\cup\{u'v',u'v'',u''v',u''v''\}$. Note that the subgraph of $G'$ induced by $u',v',u''$ and $v''$ is a 4-cycle. For convenience, we
denote $u''v''$ and $u'v'$ by $e$ and $e'$, respectively. \vskip
0.00 in

\begin{figure}
\centering
\includegraphics[width=50mm]{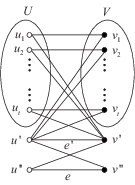}
\caption{The graph $G'$ constructed from $G$ for proving NP-completeness of conditional matching preclusion problem.} \label{g1}
\end{figure}

\vskip 0.00 in

In the following, we shall show that $G$ has a matching preclusion set of
cardinality no greater than $k$ if and only if $G'$ has a
conditional matching preclusion set of cardinality no greater than $k+1$.

\vskip 0.00 in

{\it Necessity}. Suppose that $B$ is a matching preclusion set of $G$ with $\vert B\vert\leq k$. Then $G-B$ has no perfect matchings. Let $B'=B\cup\{e\}$. We will prove that $B'$ is a conditional matching preclusion set of $G'$. Since $u'$ joins to each vertex in $V$ and $v'$ joins to each vertex in $U$, $G'-B'$ is connected. Noting that $e\in B'$, if $G'-B'$ has a perfect matching $M$, then $u'v'',u''v'\in M$. Thus, $G-B$ has a perfect matching $M\setminus\{u'v'',u''v'\}$, a contradiction. Hence, $B'$ is a conditional matching preclusion set of $G'$ with $\vert B'\vert\leq k+1$.

\vskip 0.05 in

{\it Sufficiency}. Suppose that $B'$ is a conditional matching preclusion set of
$G'$ such that $\vert B'\vert\leq k+1$. That is, $G'-B'$ has no perfect matchings, and each component of $G'-B'$ contains at least two vertices. Let $B=B'\cap E(G)$. We consider the following two cases:

\vskip 0.05 in

\noindent{\bf Case 1:} $e\in B'$. Then $\vert B\vert\leq k$. Since each component of $G'-B'$ contains at least two vertices, $u'v'',u''v'\not\in B'$. If $G-B$ has a perfect matching $M$, then
$M\cup\{u'v'',u''v'\}$ is a perfect matching of $G'-B'$, a
contradiction. Thus, $B$ is a matching preclusion set of $G$.

\vskip 0.05 in

\noindent{\bf Case 2:} $e\not\in B'$. If $B$ is a matching preclusion set of $G$ with $\vert B\vert\leq k$, we are done. If $B$ is not a matching preclusion set of $G$, then $G-B$ has a perfect matching $M$. For convenience, we denote the edges in $M$ by $u_{i}v_{i}$ with $1\leq i\leq t$. We
claim that $e'\in B'$. Otherwise, $M\cup\{e,e'\}$ is a perfect matching
of $G'-B'$, a contradiction. Also, we claim that at least one
of $v'u_{i}$ and $u'v_{i}$ is in $B'$ for each $1\leq i\leq t$. Otherwise,
$M\cup\{e\}\cup\{u'v_{i},v'u_{i}\}\setminus\{u_{i}v_{i}\}$ is a
perfect matching of $G'-B'$, a contradiction. Thus, $t+1+\vert B\vert\leq \vert
B'\vert\leq k+1$, which implies that $t\leq k$. Since $G$ is bipartite, $\Delta(G)\leq
t\leq k$. It follows that $G$ has a trivial matching preclusion
set with cardinality no greater than $k$.

We now consider the remaining case: $B$ is a matching preclusion set of $G$ such that $\vert
B\vert=k+1$. Then $B'\subseteq E(G)$. Let $M$ be a maximum matching
of $G-B$. We claim that $\vert M\vert\leq t-2$. Otherwise, $\vert M\vert=t-1$. Suppose that $u_{i}$
($1\leq i\leq t$) and $v_{j}$ ($1\leq j\leq t$) are the only two
vertices in $G$ unsaturated by $M$. Then
$M\cup\{u_{i}v',u'v_{j},e\}$ is a perfect matching of $G'-B'$, a contradiction, and the claim holds. Let $e''$ be an arbitrary edge in $B$ and
$B_{1}=B\setminus\{e''\}$. Then $\vert B_{1}\vert=k$ and
$\nu(G-B_{1})\leq \nu(G-B)+1\leq t-1$. Thus,
$B_{1}$ is a matching preclusion set of $G$ with $\vert
B_{1}\vert \leq k$. This completes the proof.
\end{proof}

\vskip 0.1 in

Now we give the decision problem of anti-Kekul\'{e} problem as follows.

\vskip 0.00 in

\textsc{Anti-Kekul\'{e} problem:}

\vskip 0.00 in

\textit{Instance:} A nonempty graph $G=(V,E)$ having a perfect matching and a positive integer $k$.

\vskip 0.00 in

\textit{Question:} Does there exist a subset $B\subseteq E$ with $\vert
B\vert\leq k$ such that $G-B$ is connected and
$G-B$ has no Kekul\'{e} structure (perfect matching)?

\begin{corollary}{\bf .}\label{NPC-ak}
Anti-Kekul\'{e} problem on bipartite graphs is NP-complete.
\end{corollary}
\begin{proof}
We adopt all the notations defined in Theorem \ref{NPC-conditional}. If we replace the condition of $G'-B'$ having no singletons by $G'-B'$ being connected, then the proof is the same as that of Theorem \ref{NPC-conditional}.
\end{proof}

\vskip 0.05 in

For a connected graph $G$ of even order, if $G$ admits an anti-Kekul\'{e} set $F$, then $G-F$ is connected and has no perfect matchings, implying that $F$ is also an $s$-restricted matching preclusion set of $G$ for any integer $s$ with $0\leq s\leq|V(G)|-1$. On the other hand, if $G$ admits an $s$-restricted matching preclusion set for small integer $s\geq 0$, an $(s+1)$-restricted matching preclusion set of $G$ does not necessarily exist. However, we have the following result.

\begin{proposition}{\bf .}\label{s-greater-half-order} Let $G$ be a nontrivial graph of even order. If $|V(G)|/2$-restricted matching preclusion set exists in $G$, then $s$-restricted matching preclusion set exists in $G$ when $0\leq s\leq |V(G)|-1$ and $mp_{s}(G)=ak(G)$ when $|V(G)|/2\leq s\leq|V(G)|-1$.
\end{proposition}
\begin{proof}
Note that after deleting a $|V(G)|/2$-restricted matching preclusion set $F$ from $G$, each component of $G-F$ contains more than $|V(G)|/2$ vertices, making $G-F$ connected. Therefore, $F$ is an anti-Kekul\'{e} set of $G$, which implies that $F$ is also an $s$-restricted matching preclusion set when $0\leq s\leq |V(G)|-1$. Moreover, we have $mp_{s}(G)=ak(G)$ when $|V(G)|/2\leq s\leq|V(G)|-1$. Hence, the statement holds.
\end{proof}

If $G$ has no anti-Kekul\'{e} set, let $s'$ be the smallest integer such that $mp_{s'}(G)=+\infty$. That is, $s'$-restricted matching preclusion set does not exist, so $mp_{s}(G)=+\infty$ for each integer $s'\leq s\leq |V(G)|-1$ and $mp_{s}(G)<+\infty$ for each nonnegative integer $s$ smaller than $s'$. Combining with the previous proposition, we have the following result.

\begin{proposition}{\bf .}\label{relation} Let $G$ be a nontrivial graph of even order. Then
\begin{equation}\label{mp-relation-formula}
mp(G)\leq mp_1(G)\leq \cdots \leq mp_{|V(G)|/2}(G)= \cdots = mp_{|V(G)|-1}(G)=ak(G).
\end{equation}
\end{proposition}
\begin{proof} If $G$ has no anti-Kekul\'{e} set, we know that $(s'-1)$-restricted matching preclusion set exists in $G$ if $s'\geq1$. Note that $s'\leq |V(G)|/2$. Otherwise, $s'-1\geq |V(G)|/2$, which implies that $|V(G)|/2$-restricted matching preclusion set exists in $G$, leading to the existence of an anti-Kekul\'{e} set in $G$, a contradiction. Accordingly, $s$-restricted matching preclusion set exists in $G$ for each integer $0\leq s\leq s'-1$. Moreover, an $s$-restricted matching preclusion set of $G$ is clearly a special $(s-1)$-restricted matching preclusion set of $G$ with $1\leq s\leq s'-1$. Thus, $mp_{s}(G)\geq mp_{s-1}(G)$. Additionally, $mp_{s}(G)=+\infty$ for each integer $s'\leq s\leq |V(G)|-1$. If $G$ has an anti-Kekul\'{e} set, by the same reason, we have $mp_{s}(G)\geq mp_{s-1}(G)$ for any integer $1\leq s\leq |V(G)|-1$. Further, by Proposition \ref{s-greater-half-order}, we arrive at the relation (\ref{mp-relation-formula}).
\end{proof}

Note that the strict inequalities in Proposition \ref{relation} can hold. It is known that $mp(G)<mp_1(G)$ for a large number of famous interconnection networks, namely hypercube \cite{Brigham,Eddie1}, folded Petersen cube \cite{Eddie00}, cube-connected cycles $CCC_n$ for $n\geq4$\cite{Li}, de Bruijn graph $UB(n)$ for $n\geq4$ \cite{Lin}, balanced hypercube \cite{Lu} and $k$-ary $n$-cube with even $k\geq4$ \cite{Wang}. In the following, we shall give an example to show that other part of strict inequalities in Proposition \ref{relation} may also hold.

\begin{figure}[h]
\centering
\includegraphics[width=110mm]{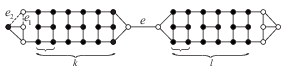}
\caption{A cubic graph with a strict inequality of (\ref{mp-relation-formula}) in which $k=3$ and $l=3$.} \label{g2}
\end{figure}

\begin{example}{\bf .} Let $G$ be a cubic graph with a cut edge $e$ shown in Fig. \ref{g2}. We note that, on the left and right sides of the cut edge, there are $k$ disjoint $P_3\Box K_2$'s and $l$ disjoint $P_3\Box K_2$'s respectively, where ``$\Box$'' means the Cartesian product of graphs, $l\geq k\geq 0$. We can check that $G$ has a perfect matching. Clearly, $G-e$ has two components of order $6k+5$ and $6l+5$ respectively, implying that it has no perfect matchings. This shows that $\{e\}$ is a minimum $s$-restricted matching preclusion set of $G$ for every $0\leq s\leq 6k+4$. On the other hand, the removal of any edge other than the cut edge from $G$ results in a connected graph with a perfect matching. Observe also that there exists two edges $e_1$ and $e_2$ (dotted lines) such that $G-\{e_1,e_2\}$ is connected and has no perfect matchings. Hence, $\{e_1,e_2\}$ is a minimum $s$-restricted matching prelusion set of $G$ for every $6k+5\leq s\leq |V(G)|-1$. So $1=mp(G)=\cdots=mp_{6k+4}(G)<mp_{6k+5}(G)=\cdots=mp_{|V(G)|/2}(G)=\cdots =mp_{|V(G)|-1}(G)=ak(G)=2$.
\end{example}

To give the complexity of $s$-restricted matching preclusion problem, we present the following decision problem.

\textsc{$s$-restricted matching preclusion problem:}

\vskip 0.00 in

\textit{Instance: A nonempty graph $G=(V,E)$ having a perfect matching, a positive integer $k$ and a positive integer $s$.}

\vskip 0.00 in

\textit{Question: Does there exist a set $B\subseteq E$ with $\vert
B\vert\leq k$ such that $G-B$ has no perfect matching and each component of $G-B$ has at least $s+1$ vertices?}

\vskip 0.05 in

\begin{theorem}{\bf .}\label{NPC-srmp}
$s$-restricted matching preclusion problem on bipartite graphs is NP-complete.
\end{theorem}
\begin{proof} We shall keep the definitions and notations introduced in Theorem \ref{NPC-conditional}.
As mentioned earlier, matching preclusion is 0-restricted matching preclusion, whose NP-completeness has already been obtained. So we may assume that $s\geq1$ in the remaining proof.
We shall prove the NP-completeness of $s$-restricted matching preclusion problem by reducing MBPMP to it in polynomial time.

In the following, we shall show that $G$ has a matching preclusion set of
cardinality no greater than $k$ if and only if $G'$ has an $s$-restricted matching preclusion set of cardinality no greater than $k+1$.

\vskip 0.00 in

{\it Necessity}. Suppose that $B$ is a matching preclusion set of $G$ with $\vert B\vert\leq k$. Then $G-B$ has no perfect matchings. Let $B'=B\cup\{e\}$. Since $u'$ joins to each vertex in $V$ and $v'$ joins to each vertex in $U$, $G'-B'$ is connected and has no perfect matchings. Thus, $B'$ is an $s$-restricted matching preclusion set of $G'$ with $\vert B'\vert\leq k+1$.

\vskip 0.05 in

{\it Sufficiency}. Suppose that $B'$ is an $s$-restricted matching preclusion set of
$G'$ such that $\vert B'\vert\leq k+1$. By Proposition \ref{relation}, $B'$ is also an $(s-1)$-restricted matching preclusion set of $G'$. Using this argument repeatedly, it follows that $B'$ is a conditional matching preclusion set of $G'$ with $\vert B'\vert\leq k+1$. By sufficiency of the proof of Theorem \ref{NPC-conditional}, it implies that there exists a  matching preclusion set $B$ of $G$ with $\vert B\vert\leq k$. This completes the proof.

\end{proof}

\section{$s$-restricted matching preclusion number}\label{sec-restricted}


In what follows, we will give some methods for determining $s$-restricted ($s\geq2$) matching preclusion numbers of regular graphs and anti-Kekul\'{e} numbers as well.




Cheng et al. \cite{Eddie1} discussed the basic obstruction to a perfect matching or an almost perfect matching in a graph with no isolated vertices. For a graph without isolated vertices, they showed that a basic obstruction to a perfect matching will be the existence of a path $uwv$, where the degree of $u$ and $v$ are 1, respectively. We define $v_{e}(G)=\min\{d_{G}(u)+d_{G}(v)-2-y_{G}(u,v)$: there exists a
vertex $w$ such that $uwv$ is a 2-path$\}$, where $d_{G}(.)$ is the
degree function and $y_{G}(u,v)=1$ if $u$ and $v$ are adjacent and
$0$ otherwise.

\begin{lemma}\cite{Eddie1}\label{v_e_G}{\bf .}
Let $G$ be a graph with an even number of vertices. Suppose every vertex in $G$ has degree at least three. Then $mp_1(G)\leq v_e(G)$.
\end{lemma}

If $mp_1(G)=v_e(G)$, then $G$ is called {\em conditionally maximally matched}. And the optimal solution of the form induced by $v_e(G)$ is called a {\em trivial} conditional matching preclusion set.

Cheng et al. \cite{Eddie0} showed that there exists a close relationship between $mp_1(G)$ and super edge-connectivity of $G$. A graph $G$ is {\em maximally edge-connected} if the edge-connectivity of $G$ is $\delta(G)$. A maximally edge-connected graph $G$ is called {\em super edge-connected} if the deletion of at most $\delta(G)$ edges results in either a connected graph or exactly two connected components, one of which is a singleton. A graph is {\em super $m$-edge-connected of order $q$} if the deletion of at most $m$ edges results in either a connected graph or a graph consisting of one big component together with a number of small components with at most $q$ vertices in total. A set of edges $F$ in a connected graph $G$ is called a {\em $g$-extra edge cut} if $G-F$ is disconnected and each remaining component of $G-F$ contains at least $g$ vertices. The {\em $g$-extra edge-connectivity} of $G$, denoted by $\lambda_g(G)$, is the minimum number of edges over all $g$-extra edge cuts of $G$. By convention, $\lambda_1(G)$ and $\lambda_2(G)$ are denoted by $\lambda(G)$ and $\lambda'(G)$, respectively. Therefore, $G$ is {\em super-$\lambda'$}, if every minimum $2$-extra edge cut isolates one edge of $G$. For $k$-regular bipartite graphs, Cheng \cite{Eddie0} presented sufficient conditions for graphs to be conditionally maximally matched as follows.

\begin{theorem}\cite{Eddie0}\label{mp1-bipartite}{\bf .} Let $G$ be a $k$-regular bipartite graph that is super ($3k-6$)-edge-connected of order 2. Then $mp_1(G)=2k-2$.
\end{theorem}

Naturally, we need the condition $k\geq3$ since this ensures the resulting graph has no isolated vertices after some edges are deleted. So we assume that $k\geq3$ in the remaining paper. Mirroring the above theorem, we obtain the following result.

\begin{theorem}\label{mps-bipartite}{\bf .} Let $G$ be a $k$-regular bipartite graph that is super ($3k-6$)-edge-connected of order 2. Then $mp_s(G)=2k-2$ for all integers $s\geq1$.
\end{theorem}

\noindent{\bf Proof.} By Theorem \ref{mp1-bipartite}, we need only to consider $s\geq2$. By Proposition \ref{relation}, we have $mp_s(G)\geq mp_1(G)=2k-2$ for all integers $s\geq2$. It remains to show that $mp_s(G)\leq 2k-2$ for all integers $s\geq2$. Let $(U,V)$ be the bipartition of $G$, where $|U|=|V|$. In addition, let $uwv$ be any 2-path in $G$. Without loss of generality, suppose $u,v\in U$ and $w\in V$. We shall show that the trivial conditional matching preclusion set obtained by $uwv$ (all the edges incident to $u$ or $v$ but not $w$), denoted by $F$, is also an $s$-restricted matching preclusion set of $G$.


If $k=3$, we shall show that $G-F$ is connected. For any $e\in F$, let $F'=F\setminus\{e\}$. Obviously, $G-F'$ contains neither isolated vertices nor isolated edges. Since $|F'|=3$ and $3k-6=3$ when $k=3$, $G-F'$ is connected. If $G-F$ is connected, we are done. So we assume that $G-F$ is disconnected, that is, $F$ is a minimal edge cut of $G$. Let $x$ be a neighbor of $w$ ($x\neq u,v$). Let $C_1$ and $C_2$ be two components of $G-F$, we may assume that $u,v,w,x\in V(C_1)$. Let $A=U\cap (V(C_1)\setminus\{u,v,x\})$ and $B=V\cap (V(C_1)\setminus\{w\})$. Since each vertex in $A$ and $B$ has degree 3, by Handshaking Lemma, we have $3|A|+d_{C_1}(u)+d_{C_1}(v)+d_{C_1}(x)=3|B|+d_{C_1}(w)$, that is, $3|A|+5=3|B|+3$, which is a contradiction. Thus, $F$ is an $s$-restricted matching preclusion set of $G$.

If $k\geq4$, the degree of each vertex (except $u$ and $v$) in $G-F$ is at least two. Then there exists no component of $G-F$ containing at most two vertices. Therefore, $G-F$ is connected since $G$ is ($3k-6$)-edge-connected of order 2 and $3k-6\geq2k-2$. Thus, $F$ is an $s$-restricted matching preclusion set of $G$. \qed

\vskip 0.1 in

We need some more definitions. An {\em independent set} in a graph $G$ is a set of vertices no two of which are adjacent. The cardinality of a maximum independent set in $G$ is called the {\em independent number} of $G$ and is denoted by $\alpha(G)$. For $X\subseteq V(G)$, we define $\gamma_{G}(X)$ as the set of edges with both ends in $X$, where the subscript $G$ will be omitted if the context is clear. Moreover, we define $\zeta(G,p,q)=\min\{\alpha(H)|H$ is an induced subgraph of $G$ with $p$ vertices and at most $q$ edges$\}$.

For non-bipartite regular graphs, Cheng et al. \cite{Eddie5} obtained the following two theorems.

\begin{theorem}\label{mp1-triangle}\cite{Eddie5}{\bf .} Let $G=(V,E)$ be a $k$-regular graph of even order, where $k\geq3$. Suppose that $G$ contains a triangle, and $G$ is $k$-edge-connected and super $(3k-8)$-edge-connected of order 2. Moreover, assume that either $|\gamma_{G}(X)|>2k-4$ for every $X\subseteq V$ of size $|X| = \frac{|V|+2}{2}$, or $\alpha(G)<\zeta(G,\frac{|V|-2}{2},2k-8)$. If $k=3$, it is additionally required that $G$ is super $(3k-7)$-edge-connected of order 2. Then $mp_1(G)=2k-3$.
\end{theorem}

\begin{theorem}\label{mp1-triangle-free}\cite{Eddie5}{\bf .} Let $G=(V,E)$ be a $k$-regular graph of even order, where $k\geq3$. Suppose that $G$ is triangle-free, and $G$ is $k$-edge-connected and super $(3k-6)$-edge-connected of order 2. Moreover, either $|\gamma_{G}(X)|>2k-3$ for every $X\subseteq V$ of size $|X| =\frac{|V|+2}{2}$, or $\alpha(G)<\zeta(G,\frac{|V|-2}{2},2k-6)$. Then $mp_1(G)=2k-2$.
\end{theorem}

We generalize the two statements above to $s$-restricted matching preclusion problem as follows.

\begin{theorem}\label{mps-triangle}{\bf .} Let $G=(V,E)$ be a $k$-regular graph of even order, where $k\geq3$. Suppose that $G$ contains a triangle, and $G$ is $k$-edge-connected and super $(3k-8)$-edge-connected of order 2. Let $s\geq2$ be any integer. If $mp_1(G)=2k-3$, then $mp_s(G)=2k-3$.
\end{theorem}

\noindent{\bf Proof.} By Proposition \ref{relation}, we have $mp_s(G)\geq 2k-3$ for all integers $s\geq2$. It remains to show that $mp_s(G)\leq 2k-3$. It suffices to find an $s$-restricted matching preclusion set $F$ with $|F|=2k-3$. Let $uwvu$ be a triangle of $G$ and let $F$ be the set of edges incident to $u$ and $v$ but not $w$. Clearly, $uv\in F$. Thus, $F$ is a trivial conditional matching preclusion set. Then $G-F$ has no perfect matchings. For convenience, let $F'=F\setminus\{uv\}$. We consider the following two cases.

\noindent{\bf Case 1.} $k=3$. Then $|F|=3$ and $|F'|=2$. Since $G$ is $3$-edge-connected, $G-F'$ is connected. Note that $uwvu$ is also a triangle in $G-F'$, by further deleting $uv$ from $G-F'$, so it is obvious that $G-F$ is connected. Therefore, $F$ is an $s$-restricted matching preclusion set.

\noindent{\bf Case 2.} $k\geq4$. Then there exist no components of $G-F'$ containing at most two vertices since the degree of each neighbor of $u$ or $v$ is at least two in $G-F'$. Since $2k-4\leq3k-8$ whenever $k\geq4$, and $G$ is super $(3k-8)$-edge-connected of order 2, $G-F'$ is connected. So $G-F$ is connected. Therefore, $F$ is an $s$-restricted matching preclusion set with $|F|=2k-3$. This completes the proof.\qed

\begin{theorem}\label{mps-triangle-free}{\bf .} Let $G=(V,E)$ be a $k$-regular graph of even order, where $k\geq3$. Suppose that $G$ is triangle-free, $k$-edge-connected and super $(3k-6)$-edge-connected of order 2. If $mp_1(G)=2k-2$, then $mp_s(G)=2k-2$ for any integer $s\geq2$.
\end{theorem}

\noindent{\bf Proof.} By Proposition \ref{relation}, we have $mp_s(G)\geq 2k-2$ for all integers $s\geq2$. It remains to show that $mp_s(G)\leq 2k-2$. It suffices to give an $s$-restricted matching preclusion set $F$ with $|F|=2k-2$. Let $uwv$ be a 2-path of $G$ and let $F$ be the set of edges incident to $u$ and $v$ but not $w$. Since $G$ is triangle-free, $uv\not\in F$. Thus, $F$ is a trivial conditional matching preclusion set. Then $G-F$ has no perfect matchings. We consider the following two cases.

\noindent{\bf Case 1.} $k=3$. Then $|F|=4$. If $G-F$ is connected, then $F$ is an $s$-restricted matching preclusion set, we are done. So we assume that $G-F$ is disconnected. Let $A=N_G(u)\cup N_G(v)\setminus\{w\}$. Then $2\leq|A|\leq 4$. We consider the following subcases.

\noindent{\bf Subcase 1.1.} $|A|=2$. We may assume that $A=\{x,y\}$. Thus, $d_{G-F}(x)=d_{G-F}(y)=1$, and $xy\not\in E(G)$. So there exists a neighbor $u'$ (resp. $v'$) of $x$ (resp. $y$) in $G-F$. We have $d_{G-F}(u')=3$ and $d_{G-F}(v')=3$. Therefore, there exist no components of $G-F$ containing at most two vertices. Let $e\in F$ be any edge and let $F'=F\setminus\{e\}$. Obviously, there exist no components of $G-F'$ containing at most two vertices. Since $G$ is super $(3k-6)$-edge-connected of order 2 and $|F'|=3$, $G-F'$ is connected. So $F$ is a minimal edge cut of $G$. Let $C_1$ and $C_2$ be two components of $G-F$, we may assume that $u,w,v\in V(C_1)$ and $x,y,u',v'\in V(C_2)$. Since $d(w)=3$, there exists a neighbor $w'$ ($w'\neq u,v$) of $w$ in $C_1$. Thus, $w$ is a cut vertex in $G$. It implies that there exist at most two edge-disjoint paths from $w$ to $x$ in $G$, which contradicts the fact that $G$ is 3-edge-connected. It follows that $G-F$ is connected.

\noindent{\bf Subcase 1.2.} $|A|=3$. We may assume that $A=\{x,y,z\}$. In addition, we assume that $xu,xv\in E(G)$. Since $G$ is triangle-free, $xy,xz\not\in E(G)$. It implies that $x$ (resp. $y$, $z$) has a neighbor of degree 3 (not $w$) in $G-F$. Then there exist no components of $G-F$ containing at most two vertices. By a similar argument of the proof of Subcase 1.1, we know that $G-F$ is connected.

\noindent{\bf Subcase 1.3.} $|A|=4$. We may assume that $A=\{x,y,x',y'\}$. In addition, we assume that $xu,yu,x'v,y'v\in E(G)$. Since $G$ is triangle-free, $xy,x'y'\not\in E(G)$. So each vertex in $A$ has degree 2 in $G-F$, indicating that there exist no components of $G-F$ containing at most two vertices. Again, we know that $G-F$ is connected.

Therefore, $F$ is an $s$-restricted matching preclusion set with $|F|=4$.

\noindent{\bf Case 2.} $k\geq4$. Noting the vertex degree of each neighbor of $u$ or $v$ is at least two in $G-F$, then there exist no components of $G-F$ containing at most two vertices. Since $2k-2\leq3k-6$ whenever $k\geq4$, combining $G$ is super $(3k-6)$-edge-connected of order 2, $G-F$ is connected. Therefore, $F$ is an $s$-restricted matching preclusion set with $|F|=2k-2$. This completes the proof.\qed

\section{Applications}

In this section, as applications of Theorems \ref{mps-bipartite}, \ref{mps-triangle} and \ref{mps-triangle-free} obtained in Section \ref{sec-restricted}, we shall determine the $s$-restricted matching preclusion numbers of complete graphs, hypercubes and hyper Petersen networks.

\subsection{Complete graphs}

\begin{theorem}\cite{Eddie1}\label{mp1-complete}{\bf .} Let $n\geq4$ be even. Then
\begin{equation*}
mp_1(K_n)=\left\{
\begin{aligned}
(n^2+2n)/8 &\ \ \   \text{if}\  n\in\{4,6,8\}, \\
2n-5 & \ \ \   \text{if}\  n\geq10.
\end{aligned}
\right.
\end{equation*}
\end{theorem}

We need the following lemma.

\begin{lemma}\cite{Wang0}\label{super-edge-connected-complete}{\bf .} Let $G$ be a complete graph with order at least four. Then $G$ is super-$\lambda'$.
\end{lemma}

Since $K_n$ ($n\geq4$) is super-$\lambda'$, we have $\lambda'(K_n)=2n-4$. Similar to Theorem \ref{mp1-complete}, we have the following result.


\begin{theorem}\label{mps-complete}{\bf .} Let $n\geq4$ be even and let $s\geq2$ be an integer. Then
\begin{equation*}
mp_s(K_n)=\left\{
\begin{aligned}
(n^2+2n)/8 &\ \ \   \text{if}\  n\in\{4,6,8\}, \\
2n-5 & \ \ \   \text{if}\  n\geq10.
\end{aligned}
\right.
\end{equation*}
\end{theorem}
\noindent{\bf Proof.} By Proposition \ref{relation}, it suffices to present an $s$-restricted matching preclusion set of size $mp_1(K_n)$. Let $F$ be a conditional matching preclusion set of $K_n$ with $|F|=mp_1(K_n)$. It follows from Theorem \ref{mp1-complete} that $|F|<2n-4$ whenever $n\geq4$. By Lemma \ref{super-edge-connected-complete}, we known that deleting less than $2n-4$ edges from $K_n$ results in a connected subgraph or a subgraph consisting of exactly two components, one of which is a singleton. It implies that $K_n-F$ is connected. Thus, $F$ is also an $s$-restricted matching preclusion set of $K_n$. \qed


\subsection{Hypercubes}

The hypercube $Q_n$ is a well-known topology for parallel computing. Any vertex $v$ of $Q_n$ is denoted by an $n$-bit binary string $v_1v_2\cdots v_n$,
where $v_i\in\{0,1\}$, for all $i$, $1\leq i\leq n$. Two vertices of $Q_n$ are adjacent if and only
if their binary strings differ in exactly one bit position.

To compute the $s$-restricted matching preclusion number of the hypercube, we need the following results.


\begin{lemma}\cite{Esfahanian}\label{lambda-2-hypercube}{\bf .} Let $n\geq2$ be a positive integer. Then $\lambda'(Q_n)=2n-2$.
\end{lemma}

\begin{lemma}\cite{Zhu}\label{lambda-3-hypercube}{\bf .} Let $n\geq2$ be a positive integer. Then $\lambda_3(Q_n)=3n-4$.
\end{lemma}

\begin{theorem}\cite{Eddie1}{\bf .} Let $n\geq3$ be an integer. Then $mp_{1}(Q_{n})=2n-2$.
\end{theorem}

\begin{theorem}{\bf .} Let $n\geq3$. For all integers $s\geq2$, $mp_{s}(Q_{n})=2n-2$.
\end{theorem}

\noindent{\bf Proof.} It is known that $Q_n$ is bipartite and $n$-regular. We only need to verify the connectivity condition in Theorem \ref{mps-bipartite}. By Lemma \ref{lambda-3-hypercube}, $\lambda_3(Q_n)=3n-4$. Let $F$ be a set of edges in $Q_n$ with $|F|\leq 3n-6$. If $Q_n-F$ is disconnected, take a largest component $C$. Then $C$ has at least three vertices; otherwise, each component of $Q_n-F$ is a singleton or doubleton. Partition the components into two classes, each having at least three vertices in total, contradicting Lemma \ref{lambda-3-hypercube}. Similarly, the other components of $Q_n-F$ other than $C$ have two vertices in total. Hence, $Q_n$ is super ($3n-6$)-edge-connected of order 2. Therefore, $mp_{s}(Q_{n})=2n-2$.  \qed


\subsection{Hyper Petersen networks}

Let $G_1=(V_1,E_1)$ and $G_2=(V_2,E_2)$ be two graphs. The {\em Cartesian product} of $G_1$ and $G_2$ is the graph $G_1\Box G_2$ whose vertex set is $V_1\times V_2$ and whose edge set is the set of all pairs $(u_1,u_2)(v_1, v_2)$ such that either $u_1v_1\in E_1$ and $u_2=v_2$, or $u_2v_2\in E_2$ and $u_1=v_1$. Das et al. \cite{Das} introduced the hyper Petersen networks $PN_n$ for $n\geq3$ as a kind of hypercube-like interconnection network. Some additional properties of $PN_n$ can be found in \cite{Ayyoub}. $PN_n$ is recursively defined as follows: $PN_3=P$, where $P$ is the Petersen graph, and $PN_n$ = $Q_{n-3}\Box P$ for $n\geq4$. Obviously, $PN_n$ is $n$-regular and has $10\times2^{n-3}$ vertices.

\begin{lemma}\cite{Eddie00}\label{mp1-Petersen}{\bf .}
$mp_1(P)=3$.
\end{lemma}

\begin{lemma}{\bf .}
If $s\geq2$, then $mp_s(P)=3$.
\end{lemma}
\noindent{\bf Proof.} Let $F$ be a conditional matching preclusion set of $P$ with $|F|=3$. It is easy to see that the Petersen graph is super edge-connected. Then $P-F$ is connected. So $F$ is also an $s$-restricted matching preclusion set. By Proposition \ref{relation}, $mp_s(P)=3$. \qed

\begin{theorem}\cite{Eddie5}\label{mp1-PNn}{\bf .}
If $n\geq4$, then $mp_1(PN_n) = 2n-2$.
\end{theorem}

\begin{theorem}{\bf .}
If $n\geq4$, then $mp_s(PN_n) = 2n-2$.
\end{theorem}
\noindent{\bf Proof.} Obviously, $PN_n$ is $n$-regular and triangle-free. By the proof of Theorem 6.2 in \cite{Eddie5}, it can be known that $PN_n$ is super edge-connected and super ($3n-6$)-edge-connected
of order 2 for $n\geq3$. By Theorems \ref{mps-triangle-free} and \ref{mp1-PNn}, we have $mp_s(PN_n) = 2n-2$.\qed

\section{Conclusions}
The MBPMP problem, arising in the structural analysis of
differential-algebraic systems, is the same as matching preclusion
problem for measuring robustness of interconnection networks. It is known that MBPMP is NP-complete, thus matching preclusion problem on bipartite graphs is also NP-complete. By reducing MBPMP to
conditional matching preclusion problem in polynomial time, we prove its NP-completeness. As a corollary, we prove NP-completeness of anti-Kekul\'{e} problem. We generalize matching preclusion and conditional matching preclusion to $s$-restricted matching preclusion and obtain its NP-completeness.

To calculate the $s$-restricted ($s\geq2$) matching preclusion numbers of graphs, we present some sufficient conditions for regular graphs and, for example, obtain $s$-restricted matching preclusion numbers for complete graphs, hypercubes, Petersen graph and hyper Petersen networks. It is interesting to study the $s$-restricted matching preclusion numbers for general graphs. Additionally, the complexity of $s$-restricted matching preclusion problem on graphs with restricted conditions, say maximum degree, planar, should be further studied.

\vskip 0.3 in

\noindent{\bf\normalsize Acknowledgments}

\vskip 0.05 in

The authors are grateful to Professor Eddie Cheng for valuable suggestions on sufficient conditions, as well as examples, of determining $s$-restricted matching preclusion number of graphs. The authors' thanks also goes to anonymous referees for their valuable comments and suggestions which improve our manuscript.

\vskip 0.3 in

\noindent{\bf\normalsize Funding} This research was supported by National Natural Science Foundation of China (Nos. 11801061, 12271228 and 12271229) and the Basic Research Project of Qinghai (No. 2021-ZJ-703).

\vskip 0.3 in

\noindent{\bf\normalsize Data Availability} Data sharing not applicable to this article as no datasets were generated or analysed during the current study.

\vskip 0.3 in

\noindent{\bf\large Declarations}

\vskip 0.3 in

\noindent{\bf\normalsize Competing Interests} The authors have not disclosed any competing interests.

\vskip 0.3 in

\noindent{\bf\normalsize Financial interests} The authors have no financial or proprietary interests in any material discussed in this article.

\end{document}